\documentclass[amstex,12pt,amssymb]{article}
\usepackage{mathtext}
\usepackage{amsmath}
\usepackage{amssymb}
\usepackage{amsxtra}
\usepackage{latexsym}
\usepackage{ifthen}
\usepackage{ulem}
\usepackage{mathbbol}
\usepackage{indentfirst}

\usepackage{color}

\makeatletter
\def\thmstyle{\it} % style of text in theorem environment
\def\@begintheorem#1#2{\it \trivlist \item[\hskip
        \labelsep{\bf #1\ #2.}]\thmstyle}
\def\@opargbegintheorem#1#2#3{\it \trivlist \item[\hskip
        \labelsep{\bf #1\ #2\ (#3).}]\thmstyle}
\makeatother

\begin{document}

\renewcommand{\baselinestretch}{1.25}
\newcommand {\beq}{\begin{equation}}
\newcommand {\eeq}{\end{equation}}
\newcommand{\red}[1]{{\color{red}#1}}

\newtheorem{theorem}{Theorem}[section]
\newtheorem{lemma}[theorem]{Lemma}
\newtheorem{proposition}[theorem]{Proposition}
\newtheorem{corollary}[theorem]{Corollary}
\newtheorem{conjecture}[theorem]{Conjecture}
\newtheorem{definition}[theorem]{Definition}
\newtheorem{remark}[theorem]{Remark}
\newtheorem{claim}[theorem]{Claim}

\title{Local large deviation principle for Wiener process with random resetting
%\footnote{The author is supported by RFBR according to the research project 18-01-00101}
 \author{\bf A. Logachov, O. Logachova,  A. Yambartsev
 }
}\maketitle

\begin{abstract}
We consider a class of Markov processes with resettings, where at random times, the Markov processes are restarted from a predetermined point or a region. These processes are frequently applied in physics, chemistry, biology, economics, and in population dynamics. In this paper we establish the local large deviation principle (LLDP) for the Wiener processes with random resettings, where the resettings occur at the arrival time of a Poisson process. Here, at each resetting time, a new resetting point is selected at random, according to a conditional distribution.  
\end{abstract}



\section{Introduction}
%\red{Here we will assume that all random
%elements are given on a stochastic basis $(\Omega,\mathfrak{F},\mathfrak{F}_t,\bf{P})$,
%$t\geq 0$.}

Random processes with resettings have recently found their
applications in various fields outside of mathematics. 
We will list some but not all applications of these processes: they are used in random search algorithms \cite{Ev1}--\cite{Lub}, in population dynamics \cite{Broc}--\cite{Pak}, and in biological and chemical models 
\cite{Ben}--\cite{Vis}. The majority of these applications used the Wiener processes
with resettings.  A Wiener process with resettings is defined as a solution to the following stochastic equation
\beq \label{1.2}
\xi(t)=w(t)-\int\limits_0^t\xi(s-)d\nu(s),
\eeq
where $w(t)$ denotes a Wiener process and $\nu(t)$ is a Poisson process with rate $\lambda$. 
The processes $w(t)$ and $\nu(t)$ are assumed to be independent.
The equation \eqref{1.2} corresponds to the case when at each Poisson arrival time, the Wiener process
restarts at the origin.

In most of the works where the processes of type \eqref{1.2} are considered,
the authors concentrate on the analysis of the corresponding stationary distributions or
additive and integral functionals. See \cite{Holl1} and references therein.
To the best of our knowledge, the paper of Meylahn et al. \cite{Mey} stands out as the only work where the large deviation principle (LDP) was established for integral functionals of the diffusion processes with resettings.

The main objective of this current paper is to establish a local large deviation principle (LLDP) 
for the trajectories of this type of processes. 
We believe that prior to this work there were no such results proved for the Wiener processes with resettings. 
Moreover, we prove the LLDP for the case where the resetting point is selected at random.

In what follows, we assume that all random elements considered here are in the probability space 
$\bigl(\Omega,\mathfrak{F}=\mathfrak{B} \cup \left (\cup_{t \geq 0}\mathfrak{F}_t\right),\bf{P}\bigr)$. 
Here, $\mathfrak{B}$ is Borel $\sigma$-algebra on $\mathbb{R}$, $\mathfrak{F}_t$ is the filtration induced by the trajectories of $\big(w(t),\nu(t)\big)$,
where $w(t)$ is the Wiener process, 
$\nu(t)$ denotes the Poisson process with rate $\lambda$, 
and the processes $w(t)$ and $\nu(t)$ are assumed to be independent.
We will examine the following stochastic equation
\beq \label{1.1}
\xi(t)=w(t)-\int\limits_0^t\zeta(\nu(s-),\xi(s-))d\nu(s),
\eeq
where a collection of independent $\mathfrak{B}$-measurable random variables
$\zeta(n,x)$, $n\in \mathbb{Z}^+:=\{0\}\cup\mathbb{N}$, $x\in \mathbb{R}$,
independent from $\big(w(t),\nu(t)\big)$.

Note that the Wiener processes with resetting satisfying (\ref{1.2})
are a special case of the processes evolving according to the equation (\ref{1.1}) with $\zeta(n,x)=x$.

In the modern literature on the LDP, various conditions on random processes are considered in order to obtain a rough exponential asymptotics for probabilities of rare events (see \cite{DZ,Feng}). In the studies where LDP have been proved for the solutions of stochastic differential equations containing an integral with respect to a Poisson process measure, a bound on the function $\zeta(n,x)$ is usually required, and is frequently given in the form of the Lipschitz condition  (see \cite{Makh}--\cite{Bud2}). In our case $\zeta(n,x)$ is the random function of $x$, which as we know, does not satisfy the Lipschitz condition, and also is unbounded.

We are interested in establishing the LLDP for positive and negative excursions of a process 
$$
\xi_T(t):=\frac{\xi(Tt)}{T}, \ t\in [0,1],
$$
where $T$ is an unbounded increasing parameter. This paper extends the ideas of our previous LLDP
result for the random walk with catastrophes \cite{LLY1}.

The trajectories of the process $\xi_T(\cdot)$ almost surely belong to the set
$\mathbb{D}[0,1]$ of c\'adl\'ag functions (i.e., right continuous and with a left limit).
For $f,g \in \mathbb{D}[0,1]$ let
$$
\rho(f,g)=\sup\limits_{t\in[0,1]}|f(t)-g(t)|.
$$

We recall the definition of LLDP.
\begin{definition}
 A family of random processes $\xi_T(\cdot)$ satisfies
the local large deviation principle (LLDP) on the set $G\subset \mathbb{D}[0,1]$ with a rate function
$I = I(f):\, \mathbb{D}[0,1] \rightarrow [0,\infty]$  and the normalizing function
$\psi(T)$ such that $\lim\limits_{T\rightarrow\infty}\psi(T) = \infty$
 if for any function $f \in G$, we have
\begin{equation}\label{opr}
\begin{aligned}
& \lim_{\varepsilon\rightarrow 0}\limsup_{T\rightarrow \infty}\frac{1}{\psi(T)}
\ln\mathbf{P}(\xi_T(\cdot)\in U_\varepsilon(f))
\\
&=\lim_{\varepsilon\rightarrow 0}\liminf_{T\rightarrow \infty}\frac{1}{\psi(T)}
\ln\mathbf{P}(\xi_T(\cdot)\in U_\varepsilon(f))=-I(f),
\end{aligned}
\end{equation}
where
$
U_\varepsilon(f):=\{g\in \mathbb{D}[0,1]: \ \rho(f,g)<\varepsilon\}.
$
\end{definition}
See \cite{BorMog1,BorMog2} for more details on the concept of LLDP.

%Recall also the definition of LDP.
%Denote the closure and interior of the set $B$ by  $[B]$ and $(B)$, respectively.

%\begin{definition} \label{d1.2}
%  A family of random processes $\xi_T(\cdot)$ satisfies
%LDP  on  metric space (m.s.) $(\mathbb{D}[0,1],\rho)$ with a rate function
%$I = I(f)\,:\, \mathbb{D}[0,1] \rightarrow [0,\infty]$ and the normalizing function
%$\psi(T)\,:\, \lim\limits_{T\rightarrow\infty}\psi(T) = \infty$,
%if, for any $c \geq 0 $ \ the set $\{ f \in \mathbb{D}[0,1]\,:\, I(f) \leq c \}$  is a compact set on m.s. $(\mathbb{D}[0,1],\rho)$
% and, for any set $B \in \mathfrak{B}_{(\mathbb{D}[0,1],\rho)}$ the following inequalities hold:
%$$
%\begin{aligned}
%\limsup_{T \rightarrow \infty} \frac{1}{\psi(T)} \ln \mathbf{P}(\xi_T(\cdot) \in B ) & \leq - I([B]),
%\\
%\liminf_{T \rightarrow \infty} \frac{1}{\psi(T)} \ln \mathbf{P}(\xi_T(\cdot) \in B ) & \geq -I((B)),
%\end{aligned}
%$$
%where $\mathfrak{B}_{(\mathbb{D}[0,1],\rho)}$ is the Borel $\sigma$-algebra  constructed by open
%cylindrical subsets of the space $\mathbb{D}[0,1]$, $I(B) = \inf\limits_{y \in B} I(y)$ for $B\in \mathfrak{B}_{(\mathbb{D}[0,1],\rho)}$,
%$I(\emptyset) = \infty$.
%\end{definition}

\bigskip
We let $\mathbf{p}_{\zeta(n,x)}(y)$ be the density of the random variable $\zeta(n,x)$. Furthermore, we assume that $\zeta(n,x)$ satisfies the following conditions.

\begin{enumerate}
\item[$\mathbf{A_0}$:] if $x=0$, then $\mathbf{P}(\zeta(n,x)= 0)=1$ for all $n\in \mathbb{Z}^+$;

\item[$\mathbf{A_+}$:] if $x>0$, then $\int_0^x\mathbf{p}_{\zeta(n,x)}(y)dy=1$
for all $n\in \mathbb{Z}^+$;

\item[$\mathbf{A_-}$:] if $x<0$, then $\int_x^0\mathbf{p}_{\zeta(n,x)}(y)dy=1$
for all $n\in \mathbb{Z}^+$;

\item[$\mathbf{B_+}$:] if $x>0$, then there exists $\Delta\geq 1$ such that
$\dfrac{1}{\Delta|x|}\leq\mathbf{p}_{\zeta(n,x)}(y)\leq \dfrac{\Delta}{|x|}$ holds true for all
$n\in \mathbb{Z}^+$, and for almost all $y\in [0,x]$;

\item[$\mathbf{B_-}$:] if $x<0$, then there exists $\Delta\geq 1$ such that
$\dfrac{1}{\Delta|x|}\leq\mathbf{p}_{\zeta(n,x)}(y)\leq \dfrac{\Delta}{|x|}$ holds true for all
$n\in \mathbb{Z}^+$, and for almost all $y\in [x, 0]$.
\end{enumerate}

Note that all of the above conditions hold in the case when $\zeta(n,x)$ is uniformly distributed in
the interval between $0$ and $x$ whenever $x \not=0$, and $\zeta(n,0)=0$. In this example, $\Delta=1$.

We use the following notations:
$\mathbb{C}_0[0,1]$ -- the set of continuous functions on the interval $[0,1]$, originating from zero;
$\mathbb{VC}_0[0,1]$ -- the subset of functions in $\mathbb{C}_0[0,1]$
with finite variation;  $\mathbb{VC}_0^M[0,1]$ -- the subset of all nondecreasing functions in $\mathbb{VC}_0[0,1]$;
$\mathbb{AC}_0^M[0,1]$ -- the subset of absolutely continuous functions in $\mathbb{VC}_0^M[0,1]$;
$\mathbb{AC}_0^+[0,1]$ ($\mathbb{AC}_0^-[0,1]$) -- the subset of absolutely continuous functions in $\mathbb{C}_0[0,1]$, taking positive (negative) values for $t \in (0,1]$;
$\text{V}_a^b (f)$ -- the total variation of a function $f$ over the interval $[a,b]$;
$\overline{B}$ -- complement of the set $B$;
$\mathbb{1}_B(\cdot)$ -- indicator function of the set  $B$;
$\lfloor a\rfloor$ -- the integer part of a number $a$.

Further in Section 2 we formulate our main results; in Section 3 we prove the LLDP; some auxiliary results are proved in Section 4.

\section{Main results}

%Here we state main theorems.

Note that if the conditions $\mathbf{A_0}$, $\mathbf{A_+}$, $\mathbf{A_-}$ hold, then 
for any $T>0$ the equation (\ref{1.1}) has a solution on the interval $[0,T]$, and it is unique.
It is also easy to prove the following result about an asymptotic upper bound for the 
maximum value of the process. 

\begin{theorem} \label{t2.2} Let the conditions $\mathbf{A_0}$, $\mathbf{A_+}$, $\mathbf{A_-}$ hold and let
the increasing function $\varphi(T)$ satisfies
$$
\lim\limits_{T\rightarrow\infty}\frac{\varphi(T)}{\sqrt{\ln(\ln T)}}=\infty.
$$
Then for any $\varepsilon>0$ %the following holds
$$
\mathbf{P}\bigg(\lim\limits_{T\rightarrow\infty}\sup\limits_{t\in[0,1]}\bigg|\frac{\xi(Tt)}{\sqrt{T}\varphi(T)}\bigg|>\varepsilon\bigg)=0.
$$
\end{theorem}
The proof of Theorem~\ref{t2.2} is quite trivial, and we omit it.

\bigskip

\noindent
To formulate our main result, we recall that any absolutely continuous function starting from zero can be uniquely represented as
a difference of functions
$f^+\in \mathbb{AC}_0^M[0,1]$ and $f^-\in \mathbb{AC}_0^M[0,1]$
such that
$$\text{V}_0^1(f) = \text{V}_0^1 (f^+) +\text{V}_0^1 (f^-).$$
The functions $f^+$ and $f^-$ are called respectively positive and negative
variations of the function $f$, see for example \cite[Ch. 1, \S 4]{Rise}.

\begin{theorem} (LLDP) \label{t2.3}
Let the conditions $\mathbf{A_0}$, $\mathbf{A_+}$, $\mathbf{A_-}$,  $\mathbf{B_+}$ hold,
then the family of the random processes $\xi_T(\cdot)$ satisfies
LLDP on the set $\mathbb{AC}_0^+[0,1]$ with normalized function
$\psi(T)=T$ and the rate function
$$
I(f)= \lambda +\frac{1}{2}\int_0^1(\dot{f}^+(t))^2dt.
$$
\end{theorem}

\begin{theorem} (LLDP) \label{t2.4}
Let the conditions $\mathbf{A_0}$, $\mathbf{A_+}$, $\mathbf{A_-}$, $\mathbf{B_-}$ hold,
then the family of the random processes $\xi_T(\cdot)$ satisfies
LLDP on the set $\mathbb{AC}_0^-[0,1]$ with normalized function
$\psi(T)=T$ and the rate function
$$
I(f)= \lambda +\frac{1}{2}\int_0^1(\dot{f}^-(t))^2dt.
$$
\end{theorem}

\noindent
Theorems~\ref{t2.3} and \ref{t2.4} implies the following form of the rate function on the set $\mathbb{AC}_0[0,1]$.
\begin{remark}
Let all the conditions $\mathbf{A_0}$, $\mathbf{A_+}$, $\mathbf{A_-}$, $\mathbf{B_+}$, $\mathbf{B_-}$ hold,
then the family of the random processes $\xi_T(\cdot)$ satisfies
LLDP on the set $\mathbb{AC}_0[0,1]$ with normalized function
$\psi(T)=T$ and the rate function
$$
I(f)=
\lambda +\frac{1}{2} \int_0^1 \bigl( \dot{f}^+(t) \mathbb{1}_{\{f(t) \ge 0\}}(t) + \dot{f}^-(t) \mathbb{1}_{\{f(t) < 0\}}(t) \bigr)^2 dt
$$
where $f(t) = 0$ at finitely  many points in $[0, 1]$.
\end{remark}

\noindent
Moreover, the proof of the theorems provides the LLDP and the corresponding rate function for the Wiener 
processes with resetting to the origin. Note that in this case the variables $\zeta(n, x)$ are deterministic functions $\zeta(n, x) = x$, and the conditions $\mathbf{B_\pm}$ do not hold. 
\begin{remark}
LLDP for positive and negative excursions of the equation $(\ref{1.2})$
(case of the deterministic resetting at zero) is a trivial task. In this case
the random process $\xi_T(\cdot)$ will stay in the neighborhood of the function
$f\in \mathbb{AC}_0^+[0,1]$ or $f\in \mathbb{AC}_0^-[0,1]$ only if the normalized
Wiener process will stay in this neighborhood, and the Poisson process will not have jumps on the interval
$[0,T]$. Thus, thanks the independence of $w(t)$ and $\nu(t)$ the rate function takes the form
$$
I(f)= \lambda +\frac{1}{2}\int_0^1(\dot{f}(t))^2dt.
$$
\end{remark}

Note that we cannot obtain the LDP for the family $\xi_T(\cdot)$ in the metric space
$(\mathbb{D}[0,1],\rho_S)$, where $\rho_S$ is Skorohod's metric. Because one can show
that the corresponding family of measures is not exponentially tight (see \cite{DZ}, Remark (a), p.8).

\section{Proof of Theorem~\ref{t2.3} and \ref{t2.4}}

\noindent
By (\ref{1.1}) the process $\xi_T(t)$ can be written as
\beq \label{4.1}
\xi_T(t)=\frac{w(Tt)}{T}-\frac{1}{T}\int\limits_0^{Tt}\zeta(\nu(s-),\xi(s-))d\nu(s):=w_{T}(t)-\xi_{T}^-(t).
\eeq
Let us first bound  $\mathbf{P}(\xi_T(\cdot)\in U_\varepsilon(f))$ from above.
For any $c>0$ and $\delta\in (0,1)$ we have
$$
\begin{aligned}
\mathbf{P} \bigl( \xi_T(\cdot)\in U_\varepsilon(f) \bigr)
\leq &\ \mathbf{P}\bigg(\sup\limits_{t\in[\delta,1]}|\xi_T(t)-f(t)|<\varepsilon,A_c\bigg)
\\ & +\mathbf{P}\bigg(\sup\limits_{t\in[\delta,1]}|\xi_T(t)-f(t)|<\varepsilon,\overline{A_c}\bigg)
:=\mathbf{P}_1+\mathbf{P}_2,
\end{aligned}
$$
where
$
A_c:=\bigl\{\omega:\nu(T)-\nu(\delta T)\leq cT\bigr\}.
$
We bound  $\mathbf{P}_1$ from above. %Using (\ref{4.1}) obtain
%$$
%\mathbf{P}_1=\mathbf{P}\bigg(\sup\limits_{t\in[\delta,1]}|w_{T}(t)-\xi_{T}^-(t)-f(t)|<\varepsilon,A_c\bigg).
%$$
Denote
$$
B_f:=\Bigl\{g\in \mathbb{VC}_0[0,1]: \dot{g}^+(t)\geq \dot{f}^+(t) \text{ for almost all } t\in[0,1]\Bigr\},
$$
where $g^+\in\mathbb{VC}_0^M[0,1]$ is positive variation of the function $g$. 
For any $r>0$ the following inequality holds
$$
\begin{aligned}
&\mathbf{P}_1=\mathbf{P}\bigg(\sup\limits_{t\in[\delta,1]}|w_{T}(t)-\xi_{T}^-(t)-f(t)|<\varepsilon,A_c\bigg)
\\
&\leq\mathbf{P}\bigg(\sup\limits_{t\in[\delta,1]}|w_{T}(t)-\xi_{T}^-(t)-f(t)|<\varepsilon,A_c,w_{T}\in K_ r^\varepsilon\bigg)+
\mathbf{P}(w_{T}\in\overline{K_ r^\varepsilon}):=\mathbf{P}_{11}+\mathbf{P}_{12},
\end{aligned}
$$
where
$$
K_ r^\varepsilon:=
\Bigl\{v\in \mathbb{C}_0[0,1]:\inf\limits_{g\in K_r}\sup\limits_{t\in[0,1]}|g(t)-v(t)|\leq
\varepsilon\Bigr\},
\ \ \ K_ r:=\big\{g:I_1(g)\leq r\big\},
$$
and the functional 
\begin{equation}\label{rateI1}
I_1(g):= \left\{
           \begin{aligned}%{rl}
                              \frac{1}{2}\int_0^1(\dot{g}(t))^2dt, & \text{ if}\ g\in \mathbb{AC}_0[0,1],\\
                              \infty,& \text{ otherwise}.
                              \end{aligned}
                              \right.
\end{equation}

Now we bound $\mathbf{P}_{11}$ from above. 
Since the random process $\xi_{T}^-(t)$ does not decrease on the interval $[\delta,1]$,
and since the set $K_r$ is a compact, from Lemma~\ref{l5.21} it follows that
there exists $\gamma(\varepsilon)>0$ such that $\gamma(\varepsilon)\rightarrow 0$
when $\varepsilon\rightarrow 0$ and
$$
\mathbf{P}_{11} \leq \mathbf{P}(w_{T}\in B_f^{\delta,\gamma(\varepsilon)} \cap K_ r^\varepsilon,A_c)\leq
\mathbf{P}(w_{T}\in B_f^{\delta,\gamma(\varepsilon)},A_c),
$$
where
$$
B_f^{\delta,\gamma(\varepsilon)}:=\bigg\{v\in \mathbb{C}_0[0,1]:\inf\limits_{g\in B_f}\sup\limits_{t\in[\delta,1]}|g(t)-v(t)|\leq\gamma(\varepsilon)\bigg\}.
$$
Thanks of independence of the processes $w(t)$ and $\nu(t)$ we obtain
$$
\mathbf{P}_{11}\leq\mathbf{P}(w_{T}\in B_f^{\delta,\gamma(\varepsilon)},A_c)=\mathbf{P}(w_{T}\in B_f^{\delta,\gamma(\varepsilon)})\mathbf{P}(A_c).
$$
Thus, for all $r>0$
\beq \label{4.2}
\mathbf{P}_{1}\leq
\mathbf{P}(w_{T}\in B_f^{\delta,\gamma(\varepsilon)})\mathbf{P}(A_c)+\mathbf{P}_{12}
=\mathbf{P}(w_{T}\in B_f^{\delta,\gamma(\varepsilon)})\mathbf{P}(A_c)+\mathbf{P}(w_{T}\in\overline{K_ r^\varepsilon}).
\eeq

%Using (\ref{4.1}), since the random process $\xi_{T}^-(t)$ does not decrease on the interval $[\delta,1]$,
%then from Lemma~\ref{l5.2} it follows that
%$$
%\mathbf{P}_1=\mathbf{P}\bigg(\sup\limits_{t\in[\delta,1]}|w_{T}(t)-\xi_{T}^-(t)-f(t)|<\varepsilon,A_c\bigg)
%\leq \mathbf{P}(w_{T}\in B_f^{\delta,\varepsilon},A_c),
%$$
%where
%$$
%B_f^{\delta,\varepsilon}:=\bigg\{v\in \mathbb{C}_0[0,1]:\inf\limits_{g\in B_f}\sup\limits_{t\in[\delta,1]}|g(t)-v(t)|\leq\varepsilon\bigg\}.
%$$
%Since $w(t)$ and $\nu(t)$ are independent, then
%\beq \label{4.2}
%\mathbf{P}_1\leq\mathbf{P} \bigl(w_{T}\in B_f^{\delta,\varepsilon},A_c\bigr)=\mathbf{P} \bigl(w_{T}\in B_f^{\delta,\varepsilon} \bigr)\mathbf{P}(A_c).
%\eeq

We bound $\mathbf{P}_2$ from above.
Denote $\tau_{k_1},\dots,\tau_{k_{\lfloor cT \rfloor}}$ the first $\lfloor cT \rfloor$ jumps of the process $\nu(Tt)$
which belong to the interval $[\delta,1]$. Denote
$$
\begin{aligned}
G_{k_l} &:=\bigl\{\omega:\xi(\tau_{k_l}-)\in[T(f(\tau_{k_l})-\varepsilon);T(f(\tau_{k_l})+\varepsilon)]\bigr\}, \ 1\leq l \leq \lfloor cT \rfloor,
\\
H_{k_l} &:=\bigl\{\omega:\zeta(k_l-1,\xi(\tau_{k_l}-))< 2 T\varepsilon\bigr\}, \ 1\leq l \leq \lfloor cT \rfloor.
\end{aligned}
$$
If a trajectory of the process $\xi_T(t)$ does not leave the set
$U_\varepsilon(f)$, then
$\zeta(k_l-1,\xi(\tau_{k_l}-))< 2 T\varepsilon$ for $\tau_{k_l}\in[\delta,1]$,
$1\leq l \leq \lfloor cT \rfloor$.
Therefore, the following inequality holds true
$$
\begin{aligned}
\mathbf{P}_2 &=\mathbf{P}\bigg(\sup\limits_{t\in[\delta,1]}|\xi_T(t)-f(t)|<\varepsilon,\overline{A_c}\bigg)
\\
&\leq\sum\limits_{r=\lfloor cT \rfloor}^{\infty}\mathbf{P}\bigg(\bigcap\limits_{l=1}^{\lfloor cT \rfloor}H_{k_l},\bigcap\limits_{l=1}^{\lfloor cT \rfloor}G_{k_l} \
\bigg| \ \nu(T)- \nu(\delta T)=r\bigg)
\mathbf{P}\bigl(\nu(T)- \nu(\delta T)=r\bigr).
\end{aligned}
$$
Denote
$
m_\delta:=\min\limits_{t\in[\delta,1]}f(t).
$
The following inequality holds 
%\begin{lemma} \label{l5.6} The inequality
\begin{equation}\label{l5.6} %$$
\mathbf{P}\bigg(\bigcap\limits_{l=1}^{\lfloor cT \rfloor}H_{k_l},\bigcap\limits_{l=1}^{\lfloor cT \rfloor}G_{k_l} \
\bigg| \ \nu(T)- \nu(\delta T)=r\bigg)\leq
\bigg(\frac{2\varepsilon\Delta}{m_\delta-\varepsilon}\bigg)^{\lfloor cT \rfloor}.
\end{equation}%$$
%holds with $H_{k_l}$,  $G_{k_l}$, $r$, $m_\delta$ defined in Section 3.
%\end{lemma}
We prove it separately in Section~\ref{aux.res}, see the subsection~\ref{proof7}. Thus,
%Using Lemma~\ref{l5.6} we obtain
$$
\mathbf{P}_2\leq
\sum\limits_{r=\lfloor cT \rfloor}^{\infty}\bigg(\frac{2\varepsilon\Delta}{m_\delta-\varepsilon}\bigg)^{\lfloor cT \rfloor}
\mathbf{P}\bigl(\nu(T)- \nu(\delta T)=r\bigr)
\leq\bigg(\frac{2\varepsilon\Delta}{m_\delta-\varepsilon}\bigg)^{\lfloor cT \rfloor}.
$$
For the sufficiently small $\varepsilon$ the inequality $m_\delta>\sqrt{\varepsilon}$ holds, therefore
\beq \label{4.3}
\mathbf{P}_2
\leq\bigg(\frac{2\varepsilon\Delta}{m_\delta-\varepsilon}\bigg)^{\lfloor cT \rfloor}\leq
\bigg(\frac{2\sqrt{\varepsilon}\Delta}{1-\sqrt{\varepsilon}}\bigg)^{\lfloor cT \rfloor}.
\eeq
From (\ref{4.3}) it follows that for any $c>0$
\beq \label{4.4}
\lim\limits_{\varepsilon\rightarrow 0}\limsup\limits_{T\rightarrow \infty}\frac{1}{T}\ln\mathbf{P}_2
\leq c \lim\limits_{\varepsilon\rightarrow 0}\ln\bigg(\frac{2\sqrt{\varepsilon}\Delta}{1-\sqrt{\varepsilon}}\bigg)=
-\infty.
\eeq
It is known (see, for example,  Theorems 2.1 and 2.2  from \cite{Freid}) that Weiner process
satisfies the LDP on the metric space $(\mathbb{D}[0,1],\rho)$, where $\rho$ is the uniform metric, with the rate function (\ref{rateI1}). It implies that for any $\varepsilon>0$
\beq \label{new3}
\lim\limits_{r\rightarrow \infty}\limsup\limits_{T\rightarrow \infty}\frac{1}{T}\ln\mathbf{P}(w_{T}\in\overline{K_ r^\varepsilon})=
-\infty.
\eeq
Thus, using (\ref{4.2}), (\ref{4.4}), (\ref{new3}) 
%Lemmas~\ref{l5.3} and \ref{l5.4}, and using
and the fact that the set $B_f^{\delta,\gamma(\varepsilon)}$ is the closed set for any
$c\in(0,1)$, $\delta\in(0,1)$ we obtain
$$
\begin{aligned}
&
\lim\limits_{\varepsilon\rightarrow 0}\limsup\limits_{T\rightarrow \infty}\frac{1}{T}\ln\mathbf{P}(\xi_T(\cdot)\in U_\varepsilon(f))\leq
\lim\limits_{\varepsilon\rightarrow 0}\limsup\limits_{T\rightarrow \infty}\frac{1}{T}\ln(\mathbf{P}_{11}+\mathbf{P}_{12}+\mathbf{P}_2)
\\
& \le \lim\limits_{\varepsilon\rightarrow 0}\limsup\limits_{T\rightarrow \infty}\frac{1}{T}\ln(3\max\{\mathbf{P}_{11},\mathbf{P}_{12},\mathbf{P}_2\})\leq
\lim\limits_{\varepsilon\rightarrow 0}\Bigl(-I_1(B_f^{\delta,\gamma(\varepsilon)})
-\lambda(1-\delta)
 +\lambda(1-\delta)c - c\ln c\Bigr)
\\
&=-I_1(B_f^{\delta})
-\lambda(1-\delta)
 +\lambda(1-\delta)c - c\ln c,
\end{aligned}
 $$
where
$$
B_f^{\delta}:=\bigg\{v\in \mathbb{C}_0[0,1]:\inf\limits_{g\in B_f}\sup\limits_{t\in[\delta,1]}|g(t)-v(t)|=0\bigg\},
$$
and in the last inequality we applied the following simple inequality
\beq %\label{5.4}
\mathbf{P}\bigl(\nu(T)-\nu(\delta T)\leq cT \bigr) \leq \exp\Bigl\{-\lambda(1-\delta)T
 +\lambda(1-\delta)cT-T c\ln c \Bigr\}.
\eeq
%Thus, using inequalities (\ref{4.2}), (\ref{4.4}), Lemmas~\ref{l5.3} and \ref{l5.4}
%and the fact that the set $B_f^{\delta,\varepsilon}$ is a closed set for any $c\in(0,1)$, $\delta>0$, then we obtain
%$$
%\begin{aligned}
%\lim\limits_{\varepsilon\rightarrow 0}\limsup\limits_{T\rightarrow \infty}\frac{1}{T}
%\ln\mathbf{P}(\xi_T(\cdot)\in U_\varepsilon(f)) &\leq
%\lim\limits_{\varepsilon\rightarrow 0}\limsup\limits_{T\rightarrow \infty}\frac{1}{T}\ln(\mathbf{P}_1+\mathbf{P}_2)
%\\
%\lim\limits_{\varepsilon\rightarrow 0}\limsup\limits_{T\rightarrow \infty}\frac{1}{T}
%\ln \bigl( 2\max\{\mathbf{P}_1,\mathbf{P}_2\} \bigr) &\leq
%\lim\limits_{\varepsilon\rightarrow 0}\Bigl(-I_1(B_f^{\delta,\varepsilon})
%-\lambda(1-\delta)
% +\lambda(1-\delta)c - c\ln c\Bigr)
%\\
%&=-I_1(B_f^{\delta})
%-\lambda(1-\delta)
% +\lambda(1-\delta)c - c\ln c,
%\end{aligned}
%$$
%where
%$$
%B_f^{\delta}:=\Bigl\{v\in \mathbb{D}[0,1]:\inf\limits_{g\in B_f}\sup\limits_{t\in[\delta,1]}|g(t)-v(t)|=0\Bigr\}.
%$$
Taking the limits $\delta\rightarrow 0$ and $c\rightarrow 0$ we obtain
$$
\lim\limits_{\varepsilon\rightarrow 0}\limsup\limits_{T\rightarrow \infty}
\frac{1}{T}\ln\mathbf{P}(\xi_T(\cdot)\in U_\varepsilon(f))\leq-I_1(B_f)-\lambda=
-I_1(f^+)-\lambda.
$$
\bigskip
\noindent
To complete the proof, we bound now $\mathbf{P}(\xi_T(\cdot)\in U_\varepsilon(f))$ from below. We have
$$
\mathbf{P}_3:=\mathbf{P}\Bigl(\sup\limits_{t\in[0,1]}|w_{T}(t)-\xi_{T}^-(t)-f(t)|<\varepsilon\Bigr)
\geq\mathbf{P}\Bigl(w_{T}(\cdot)\in U_{\frac{\varepsilon}{2}}(f^+),
\xi_{T}^-(\cdot)\in U_{\frac{\varepsilon}{2}}(f^+-f)\Bigr).
$$
Note that $f^+-f\in \mathbb{AC}_0^M[0,1]$.
If $f^+-f\equiv 0$, then
$$
\mathbf{P}\Bigl(w_{T}(\cdot)\in U_{\frac{\varepsilon}{2}}(f^+),
\xi_{T}^-(\cdot)\in U_{\frac{\varepsilon}{2}}(f^+-f)\Bigr)\geq
\mathbf{P}\Bigl(w_{T}(\cdot)\in U_{\frac{\varepsilon}{2}}(f^+),
\nu(T)=0\Bigr).
$$
Therefore, since $w(t)$ and $\nu(t)$ are independent we obtain
\beq \label{4.5}
\mathbf{P}_3\geq \mathbf{P}\Bigl(w_{T}(\cdot)\in U_{\frac{\varepsilon}{2}}(f^+)\Bigr)e^{-\lambda T}.
\eeq
Let $f^+-f\not\equiv 0$. Define
$$n(\varepsilon):=\min\bigg\{n\in \mathbb{N}:\frac{M}{n}\leq\frac{\varepsilon}{8}\bigg\},$$
where $M:=\max\limits_{t\in[0,1]}(f^+(t)-f(t))=f^+(1)-f(1)$.

Since $f^+-f$ is continuous and non-decreasing function, then
there exists a finite set of points $0=t_0<t_1<\dots<t_{n(\varepsilon)}=1$ such that
the following equalities hold true
$$
f^+(t_1)-f(t_1)=\frac{M}{n(\varepsilon)}, \ f^+(t_2)-f(t_2)=\frac{2M}{n(\varepsilon)},
\dots, \ f^+(t_{n(\varepsilon)})-f(t_{n(\varepsilon)})=M.
$$
Therefore, if the random process $\nu(Tt)$ has no jumps on $[0,t_1]$ and has only one jump
in each of intervals $[t_{k-1},t_{k}]$, $2 \leq k \leq n(\varepsilon)$, and if
random variables $\zeta(k-1,\xi(\tau_k-))$ takes values from the interval
$$
\bigg(\frac{TM}{n(\varepsilon)}-2T\varepsilon^3;\frac{TM}{n(\varepsilon)}-T\varepsilon^3\bigg),
$$
then for sufficiently small $\varepsilon$ the inequality
$$
\sup\limits_{t\in[0,1]}
\big|\xi^-(t)-(f^+(t)-f(t))\big|<\frac{\varepsilon}{2}
$$
holds. Hence, for sufficiently small $\varepsilon$ the inequality
\small\beq \label{4.6}
\mathbf{P}_3\geq\mathbf{P}\bigg(w_{T}(\cdot)\in U_{\varepsilon^3}(f^+),
\xi_{T}^-(\cdot)\in U_{\frac{\varepsilon}{2}}(f^+-f)\bigg)
\geq\mathbf{P}\bigg(w_{T}(\cdot)\in U_{\varepsilon^3}(f^+),
\bigcap\limits_{k=1}^{n(\varepsilon)} A_k,\bigcap\limits_{k=1}^{n(\varepsilon)-1} B_k\bigg),
\eeq\normalsize
holds, where
$$
\begin{aligned}
A_1&:=\{\omega:\nu(Tt_1)=0\}, \ A_k:=\{\omega:\nu(Tt_k)-\nu(Tt_{k-1})=1\},
\ 2 \leq k \leq n(\varepsilon), \\
B_k&:=\bigg\{\omega:\zeta(k-1,\xi(\tau_k-))\in
\bigg(\frac{TM}{n(\varepsilon)}-2T\varepsilon^3;\frac{TM}{n(\varepsilon)}-T\varepsilon^3\bigg)\bigg\},
 \ 1 \leq k \leq n(\varepsilon)-1.
\end{aligned}
$$
From the inequality (\ref{4.6}) it follows that
$$
\mathbf{P}_3\geq\mathbf{P}\bigg(\bigcap\limits_{k=1}^{n(\varepsilon)-1} B_k \  \bigg| \
w_{T}(\cdot)\in U_{\varepsilon^3}(f^+),
\bigcap\limits_{k=1}^{n(\varepsilon)} A_k\bigg)
\mathbf{P}\bigg(w_{T}(\cdot)\in U_{\varepsilon^3}(f^+),
\bigcap\limits_{k=1}^{n(\varepsilon)} A_k\bigg).
$$
The following inequality we prove in Section~\ref{aux.res}.
%\begin{lemma}  The inequality
\begin{equation} \label{l5.5}
\mathbf{P}\bigg(\bigcap\limits_{k=1}^{n(\varepsilon)-1} B_k \  \bigg| \
w_{T}(\cdot)\in U_{\varepsilon^3}(f^+),
\bigcap\limits_{k=1}^{n(\varepsilon)} A_k\bigg)\geq \bigg(\frac{\varepsilon^3}{\Delta(f^+(1)+\varepsilon^3)}\bigg)^{n(\varepsilon)-1}.
\end{equation}%$$
%holds with $f^+$,  $A_k$, $B_k$, $n(\varepsilon)$ defined in Section 3.
%\end{lemma}
Thanks (\ref{l5.5}) it follows that
$$
\mathbf{P}_3\geq\bigg(\frac{\varepsilon^3}{\Delta(f^+(1)+\varepsilon^3)}\bigg)^{n(\varepsilon)-1}\mathbf{P}\bigg(w_{T}(\cdot)\in U_{\varepsilon^3}(f^+),
\bigcap\limits_{k=1}^{n(\varepsilon)} A_k\bigg).
$$
Since $w_{T}(t)$ and $\nu(Tt)$ are independent, then
\beq \label{4.7}
\begin{aligned}
\mathbf{P}_3 &\geq\bigg(\frac{\varepsilon^3}{2(f^+(1)+\varepsilon^3)}\bigg)^{n(\varepsilon)-1}
\mathbf{P}\bigl(w_{T}(\cdot)\in U_{\varepsilon^3}(f^+)\bigr) \mathbf{P}\bigg(
\bigcap\limits_{k=1}^{n(\varepsilon)} A_k\bigg)
\\
&=\bigg(\frac{\varepsilon^3}{2(f^+(1)+\varepsilon^3)}\bigg)^{n(\varepsilon)-1}
\mathbf{P} \bigl( w_{T}(\cdot)\in U_{\varepsilon^3}(f^+) \bigr)
(\lambda T)^{n(\varepsilon)-1}e^{-\lambda T}
\prod\limits_{k=2}^{n(\varepsilon)}(t_k-t_{k-1}).
\end{aligned}
\eeq
Using inequalities (\ref{4.5}), (\ref{4.7}) we obtain
$$
\liminf_{T\rightarrow \infty}\frac{1}{T}
\ln\mathbf{P} \bigl( \xi_T(\cdot)\in U_\varepsilon(f)\bigr) \geq -\lambda-I_1(U_{\varepsilon^3}(f^+)).
$$
Since Weiner process satisfies LDP with rate function (\ref{rateI1}) we obtain 
$$
\lim\limits_{\varepsilon\rightarrow 0}\liminf_{T\rightarrow \infty}\frac{1}{T}
\ln\mathbf{P} \bigl( \xi_T(\cdot)\in U_\varepsilon(f) \bigr) \geq \lim\limits_{\varepsilon\rightarrow 0}
\Bigl(-\lambda-I_1(U_{\varepsilon^3}(f^+))\Bigr)=
-\lambda - I_1(f^+).%\Box
$$
The proof of Theorem~\ref{t2.4} is similar, where instead of the condition
$\mathbf{B_+}$ we work with the condition $\mathbf{B_-}$. %$\Box$

\section{Auxiliary results}\label{aux.res}

%Recall $\mathbb{VC}_0[0,1]$ is the set of continuous on interval $[0,1]$ functions started from zero
%with finite variation; and denote $\mathbb{VC}_0^M[0,1]$ the subset,  $\mathbb{VC}_0^M[0,1] \subset
%\mathbb{VC}_0[0,1]$, of non-decreasing functions.
%Denote $\mathbb{D}_0^M[0,1]$ the set of c\'adl\'ag (continuous from the right and has a limit from the left)
%functions which are non-decreasing on the interval $[0,1]$ starting from the zero.

The next technical lemma will be useful for the proof of Lemma~\ref{l5.21}.
The proof of Lemma~\ref{l5.2} is quite trivial and will be omitted.
\begin{lemma} \label{l5.2}
Let the function $f\in \mathbb{AC}_0^+[0,1]$ is represented in the form
\beq \label{n5.1}
f(t)=g_1(t)-g_2(t),
\eeq
where $g_1\in \mathbb{C}_0[0,1]$ and $g_2\in \mathbb{VC}_0^M[0,1]$.
Then the function $g_1(t)$ has the finite variation and for almost all $t\in[0,1]$ the inequality
\beq \label{n5.2}
%\liminf\limits_{\delta\rightarrow 0}\frac{g_1^+(t+\delta)-g_1^+(t)}{\delta}
\dot{g}_1^+(t)
\geq \dot{f}^+(t),
\eeq
holds true, where $g_1^+(t)$ is the positive variation of the function $g_1(t)$.
\end{lemma}

%\bigskip

Denote $(\mathbb{C}[0,1],\rho)$ the space of continuous functions on the interval
$[0,1]$ with given uniform metric $\rho$. Let
$\mathbb{D}_0^M[0,1]$  be the set of c\'adl\'ag functions
(continuous from the right and has a limit from the left) starting from the zero
which are non-decreasing on the interval $[0,1]$.

Consider the family of functions $u_T(t)$, $t\in[0,1]$, $T>0$ which can be represented
in the form $u_T(t):=\tilde{u}_T(t)-\hat{u}_T(t)$,
where $\hat{u}_T\in \mathbb{D}_0^M[0,1]$, and
$\tilde{u}_T\in \mathbb{C}_0[0,1]\cap K_{(\mathbb{C}[0,1],\rho)} $, and
$K_{(\mathbb{C}[0,1],\rho)}\subset (\mathbb{C}[0,1],\rho)$ is some compact set.

%\subsection{Proof of the Lemma~\ref{l5.21}.}
\begin{lemma} \label{l5.21}
Let for a function $f\in \mathbb{AC}_0^+[0,1]$ the following holds true
\beq \label{new2}
\lim\limits_{T\rightarrow\infty}\sup\limits_{t\in[0,1]}|u_T(t)-f(t)|=0.
\eeq
Then
$$
\lim\limits_{T\rightarrow\infty}\inf_{g\in B_f}\sup\limits_{t\in[0,1]}|\tilde{u}_T(t)-g(t)|=0.
$$
\end{lemma}

\noindent
{\it Proof.} Proof by contradiction. Suppose not. Then, there exists $\gamma>0$ such that
for any $M>0$ there exists $T>M$ and
\beq \label{new1}
\inf_{g\in B_f}\sup\limits_{t\in[0,1]}|\tilde{u}_T(t)-g(t)|\geq \gamma.
\eeq
Since the family $\tilde{u}_T$ is contained in some compact set, then, if
the inequality (\ref{new1}) holds, then there exists subsequence
$T_M$ and continuous function $\tilde{g}$ such that
$$
\lim\limits_{M\rightarrow\infty}\sup\limits_{t\in[0,1]}|\tilde{u}_{T_M}(t)-\tilde{g}(t)|=0, \ \ \
\inf_{g\in B_f}\sup\limits_{t\in[0,1]}|\tilde{g}(t)-g(t)|\geq \gamma.
$$
Therefore, from (\ref{new2}) it follows that
$$
\lim\limits_{M\rightarrow\infty}\sup\limits_{t\in[0,1]}|\hat{u}_{T_M}(t)-(\tilde{g}(t)-f(t))|=0.
$$
Wherein due to  $\hat{u}_T\in \mathbb{D}_0^M[0,1]$ the function
$
\hat{g}(t):= \tilde{g}(t) - f(t)
$
should belong to the set $\mathbb{VC}_0^M[0,1]$.
Thus, $f(t)=\tilde{g}(t)-\hat{g}(t)$, where $\tilde{g}\not\in B_f$, $\hat{g}\in\mathbb{VC}_0^M[0,1]$,
which contradicts Lemma~\ref{l5.2}.$\Box$

%\bigskip
%
%\begin{lemma} \label{l5.3}
%A family of processes $w_{T}(t)$ satisfy LDP in the metric space
%$(\mathbb{D}[0,1],\rho)$ with normalized function
%$\psi(T)=T$ and the rate function $I_1(f)$ defined in Section 3.
%\end{lemma}
%
%The proof of the Lemma \ref{l5.3} is quite trivial, so we omitted them.
%
%\bigskip
%
%\begin{lemma} \label{l5.4} For any $c\in[0,1)$, $\delta>0$ the inequality
%\beq \label{5.4}
%\mathbf{P}\bigl(\nu(T)-\nu(\delta T)\leq cT \bigr) \leq \exp\Bigl\{-\lambda(1-\delta)T
% +\lambda(1-\delta)cT-T c\ln c \Bigr\}
%\eeq
%holds true.
%\end{lemma}
%
%The proof of the Lemma \ref{l5.4} is quite trivial, so we omitted them.
%

%\bigskip

%\begin{lemma} \label{l5.5} The inequality
%$$\mathbf{P}\bigg(\bigcap\limits_{k=1}^{n(\varepsilon)-1} B_k \  \bigg| \
%w_{T}(\cdot)\in U_{\varepsilon^3}(f^+),
%\bigcap\limits_{k=1}^{n(\varepsilon)} A_k\bigg)\geq \bigg(\frac{\varepsilon^3}{\Delta(f^+(1)+\varepsilon^3)}%\bigg)^{n(\varepsilon)-1},
%$$
%holds with $f^+$,  $A_k$, $B_k$, $n(\varepsilon)$ defined in Section 3.
%\end{lemma}

%\bigskip

%\begin{lemma} \label{l5.6} The inequality
%$$\mathbf{P}\bigg(\bigcap\limits_{l=1}^{\lfloor cT \rfloor}H_{k_l},\bigcap\limits_{l=1}^{\lfloor cT \rfloor}G_{k_l} \
%\bigg| \ \nu(T)- \nu(\delta T)=r\bigg)\leq
%\bigg(\frac{2\varepsilon\Delta}{m_\delta-\varepsilon}\bigg)^{\lfloor cT \rfloor}
%$$
%holds with $H_{k_l}$,  $G_{k_l}$, $r$, $m_\delta$ defined in Section 3.
%\end{lemma}

\subsection{Proof of inequality (\ref{l5.6}). }\label{proof7}

\noindent %{\it Proof.} 
Let $G_{k_0}:=\Omega$, $H_{k_0}:=\Omega$.
We show that the inequality
\beq \label{5.7}
\mathbf{P}_l:=\mathbf{P}\bigg(H_{k_l},G_{k_l} \
\bigg| \ \nu(T)- \nu(\delta T)=r,\bigcap\limits_{d=0}^{l-1}G_{k_d},\bigcap\limits_{d=0}^{l-1}H_{k_d}\bigg)
\leq \frac{2\varepsilon\Delta}{m_\delta-\varepsilon}
\eeq
holds for $1\leq l \leq \lfloor cT \rfloor$.
We estimate from above $\mathbf{P}_l$.

We note that, by definition, a family of random variables $\zeta(k_l-1 ,m_{k_l})$, $m_{k_l}\in \mathbb{R}$ not depends on  $w(t)$ and $\nu(t)$, $\zeta(k_{l-1}-1,m_{k_{l-1}})$, $m_{k_{l-1}}\in \mathbb{R}$, $\dots$,
  $\zeta(k_1-1 ,m_{k_1})$, $m_{k_1}\in \mathbb{R}$, and hence on $\xi(\tau_{k_1}-),\dots,\xi(\tau_{k_l}-)$.
Therefore, the next inequality
$$
\mathbf{P}\bigg(H_{k_l},G_{k_l} \
\bigg| \ \nu(T)- \nu(\delta T)=r,\bigcap\limits_{d=0}^{l-1}G_{k_d},\bigcap\limits_{d=0}^{l-1}H_{k_d}\bigg)
\leq\int\limits_0^{2T\varepsilon}
\bigg(\int\limits_{T(m_\delta-\varepsilon)}^{T(M_1+\varepsilon)}\mathbf{p}_{\zeta(k_l,x)}(y)d\tilde{F}(x)\bigg)dy
$$
holds, where $M_1:=\max\limits_{t\in[0,1]}f(t)$,
$$
\tilde{F}(x):=\mathbf{P}\bigg(\xi(\tau_{k_l}-)<x \
\bigg| \ \nu(T)- \nu(\delta T)=r,\bigcap\limits_{d=0}^{l}G_{k_d},\bigcap\limits_{d=0}^{l-1}H_{k_d}\bigg).
$$
Using the condition $\mathbf{B}_+$, we get for sufficiently small $\varepsilon$
$$
\begin{aligned}
\mathbf{P}_l & \leq \int\limits_0^{2T\varepsilon}
\bigg(\int\limits_{T(m_\delta-\varepsilon)}^{T(M_1+\varepsilon)}
\mathbf{p}_{\zeta(k_l-1,x)}(y)d\tilde{F}(x)\bigg)dy\leq
\int\limits_0^{2T\varepsilon}
\bigg(\int\limits_{T(m_\delta-\varepsilon)}^{T(M_1+\varepsilon)}
\frac{\Delta}{|x|}d\tilde{F}(x)\bigg)dy
\\
&\leq
\int\limits_0^{2T\varepsilon}
\bigg(\int\limits_{T(m_\delta-\varepsilon)}^{T(M_1+\varepsilon)}
\frac{\Delta}{T(m_\delta-\varepsilon)}d\tilde{F}(x)\bigg)dy\leq \frac{2T\varepsilon\Delta}{T(m_\delta-\varepsilon)}
=\frac{2\varepsilon\Delta}{m_\delta-\varepsilon}.
\end{aligned}
$$
Thus, the inequality (\ref{5.7}) is proved. Using the inequality
(\ref{5.7}), we obtain
$$
\mathbf{P}\bigg(\bigcap\limits_{l=1}^{\lfloor cT \rfloor}H_{k_l},\bigcap\limits_{l=1}^{\lfloor cT \rfloor}G_{k_l} \
\bigg| \ \nu(T)- \nu(\delta T)=r\bigg)=
\prod\limits_{l=1}^{\lfloor cT \rfloor}\mathbf{P}_l\leq
\bigg(\frac{2\varepsilon\Delta}{m_\delta-\varepsilon}\bigg)^{\lfloor cT \rfloor}.%\Box
$$

\subsection{Proof of inequality (\ref{l5.5})}

\noindent
We show that, for $1\leq k \leq n(\varepsilon)-1$ the inequality
\beq \label{5.5}
\mathbf{P}_k:=\mathbf{P}\bigg( B_k \  \bigg| \
w_{T}(\cdot)\in U_{\varepsilon^3}(f^+),
\bigcap\limits_{r=1}^{n(\varepsilon)} A_r,B_1,\dots,B_{k-1}\bigg)\geq\frac{\varepsilon^3}{\Delta(f^+(1)+\varepsilon^3)}
\eeq
holds. If events $\{\omega:w_{T}(\cdot)\in U_{\varepsilon^3}(f^+)\}$,
$\bigcap\limits_{r=1}^{n(\varepsilon)} A_r$, $B_1,\dots,B_{k-1}$ have occurred, then
\beq\label{5.6}
\begin{aligned}
&T (f^+(1)+\varepsilon^3)>T(f^+(\tau_{k})+\varepsilon^3)>\xi(\tau_{k}-)\geq T\bigg(w_T(\tau_{k}-)-(k-1)\bigg(\frac{M}{n(\varepsilon)}-\varepsilon^3\bigg)\bigg)
\\
&>T\bigg(f^+(\tau_{k})-\varepsilon^3-(k-1)\bigg(\frac{M}{n(\varepsilon)}-\varepsilon^3\bigg)\bigg)
>T\bigg(f^+(t_{k})-\varepsilon^3-(k-1)\bigg(\frac{M}{n(\varepsilon)}-\varepsilon^3\bigg)\bigg)
\\
& >T\bigg(f^+(t_{k})-f(t_{k})-\varepsilon^3-(k-1)\bigg(\frac{M}{n(\varepsilon)}-\varepsilon^3\bigg)\bigg)
>\frac{TM}{n(\varepsilon)}-T\varepsilon^3.
\end{aligned}
\eeq
We note that, by definition, the family of
random variables $\zeta(k-1 ,m_k)$, $m_k\in \mathbb{R}$
not depends on $w(t)$ and $\nu(t)$, $\zeta(k-2 ,m_{k-1})$, $m_{k-1}\in \mathbb{R}$, $\dots$,
  $\zeta(0,m_1)$, $m_1\in \mathbb{R}$, and hence on $\xi(\tau_k-),\dots,\xi(\tau_1-)$.
Therefore, using inequality (\ref{5.6}), we obtain
$$
\begin{aligned}
\mathbf{P}_k & =\mathbf{P}\bigg( B_k \  \bigg| \
w_{T}(\cdot)\in U_{\varepsilon^3}(f^+),
\bigcap\limits_{r=1}^{n(\varepsilon)} A_r,B_1,\dots,B_{k-1}\bigg)
\\
&=\mathbf{P}\bigg( \zeta(k-1,\xi(\tau_k-))\in
\bigg(\frac{TM}{n(\varepsilon)}-2T\varepsilon^3;\frac{TM}{n(\varepsilon)}-T\varepsilon^3\bigg) \  \bigg| \
w_{T}(\cdot)\in U_{\varepsilon^3}(f^+),
\bigcap\limits_{r=1}^{n(\varepsilon)} A_r,\bigcap\limits_{r=1}^{k-1} B_r\bigg)
\\
&=\int\limits_{\frac{TM}{n(\varepsilon)}-2T\varepsilon^3}^{\frac{TM}{n(\varepsilon)}-T\varepsilon^3}
\bigg(\int\limits_{\frac{TM}{n(\varepsilon)}-T\varepsilon^3}^{T(f^+(1)+\varepsilon^3)}\mathbf{p}_{\zeta(k-1,x)}(y)dF(x)\bigg)dy,
\end{aligned}
$$
where
$$
F(x):=\mathbf{P}\bigg(\xi(\tau_{k}-)<x \
\bigg| \ w_{T}(\cdot)\in U_{\varepsilon^3}(f^+),
\bigcap\limits_{r=1}^{n(\varepsilon)} A_r,\bigcap\limits_{r=1}^{k-1} B_r\bigg).
$$
Using the condition $\mathbf{B}_+$, we get for sufficiently large $T$
$$
\begin{aligned}
\mathbf{P}_k & \geq \int\limits_{\frac{TM}{n(\varepsilon)}-2T\varepsilon^3}^{\frac{TM}{n(\varepsilon)}-T\varepsilon^3}
\bigg(\int\limits_{\frac{TM}{n(\varepsilon)}-T\varepsilon^3}^{T(f^+(1)+\varepsilon^3)}\frac{1}{\Delta|x|}dF(x)\bigg)dy
\geq \int\limits_{\frac{TM}{n(\varepsilon)}-2T\varepsilon^3}^{\frac{TM}{n(\varepsilon)}-T\varepsilon^3}
\frac{1}{\Delta T(f^+(1)+\varepsilon^3)}dy
\\
& =\frac{T\varepsilon^3}{\Delta T(f^+(1)+\varepsilon^3)}= \frac{\varepsilon^3}{\Delta (f^+(1)+\varepsilon^3)}.
\end{aligned}
$$
Thus, the inequality (\ref{5.5}) is proved. Using the inequality
(\ref{5.5}), we get
$$
\mathbf{P}\bigg(\bigcap\limits_{k=1}^{n(\varepsilon)-1} B_k \  \bigg| \
w_{T}(\cdot)\in U_{\varepsilon^3}(f^+),
\bigcap\limits_{k=1}^{n(\varepsilon)} A_k\bigg)=
\prod\limits_{k=1}^{n(\varepsilon)-1}\mathbf{P}_k\geq
\bigg(\frac{\varepsilon^3}{\Delta(f^+(1)+\varepsilon^3)}\bigg)^{n(\varepsilon)-1}.%\Box
$$

\section*{Acknowledgments}

The authors would like to thank the anonymous referee for providing valuable comments and feedback that helped us 
greatly in improving the manuscript.

This work is supported by FAPESP grant 2017/20482-0.
LA thanks RSF project 18-11-00129 and Institute of Mathematics and
Statistics of University of S\~ao Paulo for hospitality. AY thanks CNPq and FAPESP for the
financial support via grants 301050/2016-3 and 2017/10555-0, respectively.

\end{document}